\documentclass[12pt]{article}
\usepackage{amsmath}
\usepackage{amssymb}
\usepackage{amsfonts, amsthm}
\begin{document}
\title{Near identity transformations for the Navier-Stokes
equations,  (version 12/01).}
\author{Peter Constantin
\\Department  of Mathematics\\The University of Chicago }
\maketitle
\newtheorem{thm}{Theorem}
\newtheorem{prop}{Proposition}

\section{Introduction}
Ordinary incompressible Newtonian fluids are described by the 
Navier-Stokes equations. These equations have been used by 
engineers and physicists with a great deal of success and the range of 
their validity and applicability is well established. Together with 
other fundamental systems like the Schr\"{o}dinger and Maxwell 
equations, these equations are among the most important equations of  
mathematical physics. Nevertheless, their mathematical theory is incomplete
and requires cut-offs. The present state of knowledge is 
such that different approximations seem to be useful for 
different purposes.  The mathematical questions of existence and 
regularity for incompressible fluid equations have been discussed in 
many books and review articles (for instance \cite{a}, \cite{cfbook}, \cite{pll}, \cite{tbook}, \cite{cnot}, \cite{c2k}, \cite{Ma}, \cite{serr}, \cite{t} and many more). In this work I describe some results reflecting 
research concerning diffusive-Lagrangian aspects of the Navier-Stokes 
equations \cite{cel}, \cite{xlns}.   
There are two distinct classes of approximations of the Navier-Stokes
equations that we consider. In one class the energy dissipation is treated 
exactly, but the vorticity equation is not exact. This class contains the
Galerkin approximations \cite{cfbook}, \cite{tbook} and mollified equations 
\cite{ckn}, \cite{pll} (see (\ref{nsmol}) below).
The other class treats the vorticity equation exactly but the energy dissipation is
approximated. This is the class of vortex methods \cite{chorin} 
(see (\ref{bob}) below) and their generalizations. This class is related by a change of variables to a class of filtered approximations (\ref{phil}) of the  formulation \cite{kuz}, \cite{ose}; the models  {\cite{mar}, \cite{holm} are a subclass of these.  The Navier-Stokes equations and their various 
approximations can be described in terms of near identity transformations. 
These are diffusive particle path transformations of physical space 
that start from the identity. The active velocity is obtained from the 
diffusive path transformation and a virtual velocity using the Weber 
formula. The active vorticity is computed from the diffusive path 
transformation and a virtual vorticity using a Cauchy formula. The path 
transformation and the virtual fields are computed in 
Eulerian coordinates (``laboratory frame''). In the absence of kinematic 
viscosity, both the virtual velocity and the virtual vorticity are 
passively transported (``frozen in'') by the flow. In the presence of 
viscosity, these fields obey diffusion equations with coefficients that are 
proportional to the kinematic viscosity and are derived from the 
diffusive transformations. The diffusive path transformations are used 
for short time intervals, as long as the transformations do not stray too much from the identity. The duration of these intervals is determined by 
the requirement of invertibility of the gradient map. If and when the 
viscosity-induced change in the Jacobian reaches a pre-assigned level, one 
stops, and one restarts the calculation from the identity transformation, 
using as initial virtual field the previously computed active field. 

\section{Energy Dissipation}
An incompressible fluid of constant density and temperature can be described in terms of the fluid velocity $u(x,t)$ and pressure $p(x,t)$, functions of 
Eulerian (laboratory) coordinates  $x$,  representing position and $t$, 
representing time. The Navier-Stokes equations are an expression
of Newton's second law, and in the absence of external sources of energy
they are
\begin{equation}
\partial_t u + u\cdot\nabla u -\nu\Delta u +\nabla p = 0, \label{nse}
\end{equation}
coupled with the constraint of incompressibility
\begin{equation}
\nabla\cdot u = 0. \label{divu0}
\end{equation}
In  the initial value problem,  a velocity
\begin{equation}
u(x,0) = u_0\label{id}
\end{equation}
is given at $t=0$. The coefficient $\nu>0$ is the kinematic viscosity.  
In an idealized situation the velocity is defined for all $x\in{\mathbf R}^3$ 
and vanishes at infinity. The vorticity equation is obtained by taking the curl of the Navier-Stokes equations:
\begin{equation}
\partial_t \omega  + u\cdot\nabla \omega -\nu\Delta \omega = \omega\cdot\nabla u
\label{omegu}
\end{equation}
where 
\begin{equation}
\omega  = \nabla\times u
\label{oo}
\end{equation}
is the vorticity. 
If $\omega$ is divergence-free one may invert this relation:
Defining a stream-vector $\psi$ that satisfies $-\Delta \psi = \omega$ using the 
Newtonian potential and then taking its curl, one obtains the familiar Biot-Savart law
\begin{equation}
u(x) = -\frac {1}{4\pi}\int _{{\mathbf R}^3}\frac{1}{|x-y|^3}(x-y)\times \omega (y) dy.
\label{bs}
\end{equation}

If the initial velocity vanishes then $u(x,t)=0$, $p(x,t)=0$ solve the
equations. Moreover, if $u_0$ is close to $u_0 = 0$ then 
the solution $u(x,t)$, $p(x,t)$ exists for all time, is smooth and converges to $0$. 
The open question in this situation concerns the behavior of the solution for 
large initial data. The Navier-Stokes equations are a parabolic regularization of the 
Euler equations (obtained by setting $\nu=0$). Although the viscous term is important, for the study of large data one needs to consider 
properties of the Euler equations.  The Navier-Stokes equations conserve momentum (integral of velocity in the present setting). The total kinetic energy
\begin{equation}
\frac{1}{2}\int_{{\mathbf{R}}^3}|u(x,t)|^2dx = K(t)
\label{kt}
\end{equation} 
is dissipated by viscosity
\begin{equation}
\frac{1}{2}\int_{{\mathbf{R}}^3}|u(x,t)|^2dx + \nu\int_{t_0}^t\int_{{\mathbf R}^3}
|\nabla u(x,s)|^2dx ds \le \frac{1}{2}\int_{{\mathbf{R}}^3}|u(x, t_0)|^2dx.
\label{disspa}
\end{equation}
The dissipation of kinetic energy is the strongest source of quantitative 
information about the Navier-Stokes equations that is presently known for 
all solutions.   
This dissipation is used to construct Leray weak solutions with finite kinetic 
energy that exist for all time, $u\in L^{\infty}(dt; L^2(dx))$, $\nabla u \in
L^2(dt\otimes dx)$ \cite{l}. This class of solutions is very wide. The
solutions have partial regularity \cite{ckn} but are not known to be smooth.
The uniqueness of the Leray weak solutions is not known. Neveretheless, 
immediately after inception, at positive times arbitrarily close to the 
initial time, the Leray solutions have square integrable gradients and 
the solutions become smooth for an interval of time. The 
solution is then uniquely determined and remains smooth for an interval 
of time whose duration is bounded below by a non-zero constant. The issue 
is whether the smooth behavior continues for all time. The simplest self-similar blow up ansatz of Leray has been 
ruled out \cite{ne}, \cite{ts}. The most 
important task is to obtain good a priori bounds for smooth solutions 
of the Navier-Stokes equations. If one has good bounds then 
the smoothness and uniqueness of the solution can be shown to 
persist \cite{serr}, \cite{serrin}, \cite{tbook}, \cite{cfbook}, \cite{pll}.
 In situations in which such bounds are not  available, the study of solutions of the 
Navier-Stokes equations needs to be pursued by considering long-lived approximate solutions. The advantage of dealing with approximations, besides practicality,  is conceptual simplicity: one may formulate sufficient conditions for global regularity quantitatively, in terms of the approximate solutions. If one devises approximations and obtains uniform bounds for  them, then, by removing the approximation, one obtains rigorous bounds for weak solutions of the Navier-Stokes equations that are valid for all time. For instance one can prove:

\begin{thm} Let $u_0$ be a function in $L^2({\mathbf R}^3)$, that
satifies the divergence-free condition $\nabla\cdot u_0 = 0$ in the sense
of distributions. Let $T>0$ be arbitrary. There exists a Leray weak 
solution $(u(x,t), (p(x,t))$ of the Navier-Stokes equations that 
is defined for $t\in [0,T]$, satisfies
$$
\frac{1}{2}\int_{\mathbf R^3} |u(x,t)|^2 dx  + \nu \int_{t_0}^t\int_{{\mathbf R}^3}|\nabla u (x,s)|^2 dxds \le \frac{1}{2}\int_{\mathbf R^3} |u(x,t_0)|^2 dx
$$
for all $t\ge t_0$ and $t_0\in I \subset [0, T]$, where $I$ is a set of full measure $|I | = T$. The initial time belongs to it, $0\in I$, in other words $t_0 = 0$ is allowed.  
In addition, the solution satisfies
$$
\int_0^T \|u(\cdot, s)\|_{L^{\infty}}ds \le K_{\infty}
$$
with $K_{\infty}$ a length scale determined mainly by the initial kinetic
energy and viscosity:
$$
K_{\infty} = C\left\{ \nu^{-2}\int_{{\mathbf R}^3}|u_0(x)|^2dx + \sqrt{\nu T}\right\}.
$$
If the initial vorticity $\omega_0 = \nabla\times u_0$ is in $L^1$, 
$\int_{{\mathbf R}^3}|\omega_0|dx <\infty$ then it remains bounded in 
$L^1$ and moreover
$$
\int_{\mathbf R^3}|\omega (x,t)|dx + \frac{1}{2\nu}\int_{\mathbf R^3} |u(x,t)|^2dx \le
\int_{\mathbf R^3}|\omega(x,t_0)|dx + \frac{1}{2\nu}\int_{\mathbf R^3} |u(x,t_0)|^2dx
$$
for $t\ge t_0$, $t_0\in I$. If the initial data is in $H^1$ and the initial 
Reynolds number is small then the solution is infinitely differentiable for 
positive time and converges to $0$. More precisely, if
$$
R_0 = \frac{1}{\nu}\left(\int_{{\mathbf R}^3}|u_0(x)|^2dx\right )^{\frac{1}{4}}
\left (\int_{{\mathbf R}^3}|\nabla u_0(x)|^2dx\right )^{\frac{1}{4}} <
2\pi^{-\frac{1}{2}}3^{-\frac{3}{4}}
$$
then the solution exists for all $t>0$,  belongs to
${\mathcal C}^{\infty}({\mathbf R}^3)$ and converges to
$0$.

\end{thm}

\vspace{.5cm}

This theorem combines the bound in \cite{c91} that was proved using 
a version of the retarded mollification approximation procedure of 
\cite{ckn} with the result of \cite{gft} that was proved
in the space-periodic case using Galerkin approximations. The last statement
about smooth solutions is proved by studying the evolution of the
product of energy and enstrophy. The specific constant (about .495) 
comes from the fact that, for divergence-free, zero-mean functions $\|u\|_{L^{\infty}} \le
\sqrt{4\pi}\sqrt{\|\nabla u\|_{L^2}}\sqrt{\|\Delta u\|_{L^2}}$ (see
(\ref{interpol})) below. If the initial Reynolds number is small then it stays
small and its rate of dissipation is a well-known quantity that controls
global existence (see (\ref{suf1}) below). Both the 
mollification approximation and the Galerkin truncation approximation 
procedure respect the energy dissipation inequality (\ref{disspa}) 
exactly, but they introduce errors in the vorticity equation. The 
(not retarded) mollification equation is described below. 
One defines a mollified  $u$ by 
\begin{equation}
[u]_{\delta} = \int_{{\mathbf R}^3} \delta^{-3}J\left(\frac{x-y}{\delta}\right )u(y)dy   = J_{\delta}(-i\nabla )u .\label{[u]}
\end{equation}
Here $\delta >0$ and the positive kernel $J$ is 
normalized $\int_{{\mathbf R}^3}J(x)dx =1$, smooth and decays 
sufficiently fast at infinity. Two canonical examples of such $J$ are
the Poisson kernel $J(x) = \pi^{-2}(1+ |x|^2)^{-2}$ and the Gaussian $J(x) = (2\pi)^{-{3/2}}e^{-|x|^{2}/2}$. The Fourier transforms of $J$,  $\widehat{J}(\xi) = e^{-|\xi |}$ and, respectively $\widehat{J}(\xi) = exp{\left(-\frac{|\xi |^2}{2}\right )}$, are non-negative, vanish at the origin, decay rapidly and are bounded above by $1$. Because of the fact that at the Fourier transform level
one has 
\begin{equation}
\widehat{[u]_{\delta}} (\xi) = \widehat{J}(\delta \xi)\widehat{u}(\xi), \label{hat}
\end{equation}
the operator of convolution with $J_{\delta}$, $[u]_{\delta} = J_{\delta}(-i\nabla)u$ is a classical  smoothing approximation of the identity. 
The mollified equation is
\begin{equation}
\partial_t u + [u]\cdot \nabla u -\nu\Delta u + \nabla p = 0\label{nsmol}
\end{equation}
together with $\nabla\cdot u = 0$. Here $[u] =[u]_{\delta}$ is computed by applying the mollifier at each instance of time. This nonlinear partial differential equation has  global solutions  for arbitrary divergence-free initial data $u_0\in L^2$. The solutions are 
smooth on $(0,T]\times {\mathbf R}^3$ and, moreover the energy inequality 
(\ref{disspa}) is valid for any $t_0\in [0,T]$, $t\ge t_0$. The  vorticity of the mollified
equation does not obey exactly (\ref{omegu}). By contrast, classical vortex methods \cite{chorin} respect the structure of the vorticity equation (\ref{omegu}) but do not obey exactly the energy dissipation inequality (\ref{disspa}). In this paper we call vortex methods the equations 
\begin{equation}
\partial _t\omega + [u]\cdot \omega - \nu \Delta \omega  = \omega \cdot\nabla [u]\label{bob}
\end{equation}
with $u$ calculated from $\omega $ using the Biot-Savart law (\ref{bs}), and $[u] = [u]_{\delta}$ computed from $u$ using the mollifier (\ref{[u]}).  Both equation and solutions depend on $\delta$ but we will keep
notation light by dropping the reference to this dependence:  
$\omega = \omega_{\delta}, \,\, [u] = [u_{\delta}]_{\delta}$. 
These vortex methods may also be described by using
an auxilliary variable $w$. One considers the equation
\begin{equation}
\partial_t w + [u]\cdot\nabla w - \nu\Delta w + \left (\nabla [u]\right )^{*}w = 0.
\label{phil}
\end{equation}
($M^*$ means the transposed matrix). A direct calculation verifies that the curl of  $w$ , $\nabla\times w$ obeys the equation (\ref{bob}), as does $\omega$. This calculation uses only the fact that $[u]$ 
is divergence-free. The system  formed by the equation (\ref{phil}) coupled with
\begin{equation}
[u] = J_{\delta}(-i\nabla){\mathbf P}(w)\label{wudelta}
\end{equation}
is equivalent to (\ref{bob},  \ref{bs}, \ref{[u]}). Here 
${\mathbf P}$,
\begin{equation}
{\mathbf P}_{jl} = \delta_{jl} -
\partial_j\Delta^{-1}\partial_l\label{pcompo}
\end{equation}
is the Leray-Hodge projector on divergence-free vectors. The initial $w$ is
required to satisfy ${\mathbf P}w_0 = u_0$. At fixed positive $\delta$ the solution is smooth and global. If $\nu =0$ then these systems have a
Kelvin circulation theorem: the integral $\oint_{\gamma}w\cdot dx$ is
conserved along closed paths $\gamma$ that are transported by the flow of $[u]$. (In contrast, the mollified equations do not have a Kelvin
circulation theorem). The energy dissipation
principle for the vortex method is
$$
\frac{1}{2}\int_{{\mathbf R}^3}u(x,t)\cdot [u](x,t)dx + \nu \int_{t_0}^t\int_{{\mathbf R}^3}{\mathbf {tr}}\left \{(\nabla u)(x,s)(\nabla [u](x,s))^*\right\}dxds \le 
$$
\begin{equation}
\frac{1}{2}\int_{{\mathbf R}^3}u(x,t_0)\cdot [u(x,t_0)]dx\label{dissp[u]}
\end{equation} 
This is obtained by taking the scalar product of (\ref{bob}) with $[\psi]$
where $u = \nabla\times \psi$. One can obtain the energy dissipation 
principle also
by taking the scalar product of (\ref{phil}) with $[u]$. One uses the fact 
that $J_{\delta}(-i\nabla)$ is a scalar operator (multiple of the identity as a matrix, i.e. acts separately on each component of a vector) that commutes with 
differentiation. Then the cancellation of the nonlinearity follows from the divergence free condition. The energy dissipation
principle gives strong control on the mollified quantities (or weak control 
on the unmollified ones):
\begin{equation}
\frac{1}{2}\int_{{\mathbf R}^3}u(x,t)\cdot [u](x,t)dx = \frac{1}{2}\int_{{\mathbf R}^3}\left (\widehat{J}(\delta\xi)\right )^{-1} \left | \widehat{[u]}(\xi,t)\right |^2d\xi\label{u[u]}
\end{equation}
and
\begin{equation}
 \int_{{\mathbf R}^3}{\mathbf {tr}}\left \{(\nabla u)(x,s)(\nabla [u](x,s))^*\right\}dx = \int_{{\mathbf R}^3}\left (\widehat{J}(\delta\xi)\right )^{-1} |\xi|^2 \left | \widehat{[u]}(\xi,t)\right |^2d\xi\label{gradu[u]}
\end{equation}
Because ${\widehat{J}}^{-1}$ is a positive function that grows exponentially 
at infinity, the inequality implies real analytic control on $[u]$:
$$
\frac{1}{2}\int_{{\mathbf R}^3}\left |(J_{\delta}(-i\nabla))^{-\frac{1}{2}}([u])(x,t)\right |^2dx + 
$$
\begin{equation}
\int_{0}^t \int_{{\mathbf R}^3}\left |\nabla(J_{\delta}(-i\nabla))^{-\frac{1}{2}}([u])(x,s)\right |^2dxds \le \frac{1}{2}\int_{{\mathbf R}^3}|u_0(x)|^2dx\label{en[u]}
\end{equation}
One needs to bear in mind, however, that this is a weaker
bound than the bound provided by the energy dissipation (\ref{disspa}) 
for the mollified equation (\ref{nsmol}), where $\left (J_{\delta}(-i\nabla)\right )^{-1}[u]$ is bounded in $L^2$.

\section{Uniform bounds}
The energy dissipation principle (\ref{disspa}) holds exactly for the mollified equation (\ref{nsmol}) and has a counterpart for the 
vortex method (\ref{bob}) in (\ref{dissp[u]}, \ref{en[u]}). These are uniform 
inequalities, in the sense that the coefficients are $\delta$-independent 
and the right hand sides are bounded uniformly for all $\delta>0$. Most 
uniform bounds are inherited by the solution of the Navier-Stokes equations 
by passage to limit. Some
uniform bounds for the equation (\ref{nsmol}) can be summarized as follows:

\begin{thm} Let $u_0$ be a square-integrable, divergence-free function.
Let $\delta >0$. Then there exits a unique
solution $(u,p)$ of (\ref{nsmol}) defined for all $t>0$. The solution is 
real analytic for positive times. The limit $\lim_{t\to 0}u(x,t) = u_0(x)$
holds in a weak sense in $L^2$. The energy inequality (\ref{disspa}) holds for
any $0\le t_0\le t$. The uniform bound
$$
\int_0^T \|u(\cdot, s)\|_{L^{\infty}}ds \le K_{\infty}
$$
holds with
$$
K_{\infty} = C\left\{ \nu^{-2}\int_{{\mathbf R}^3}|u_0(x)|^2dx + \sqrt{\nu T}\right\}
$$
and $C$ a universal constant, independent of $\delta$.
If the initial vorticity $\omega_0 = \nabla\times u_0$ is in $L^1$, 
$\int_{{\mathbf R}^3}|\omega_0|dx <\infty$ then it remains bounded in 
$L^1$ and moreover
$$
\int_{\mathbf R^3}|\omega (x,t)|dx + \frac{1}{2\nu}\int_{\mathbf R^3} |u(x,t)|^2dx \le
\int_{\mathbf R^3}|\omega(x,t_0)|dx + \frac{1}{2\nu}\int_{\mathbf R^3} |u(x,t_0)|^2dx
$$
for all $t\ge t_0$. In addition the vorticity direction
$$
\xi (x,t) = \frac{\omega(x,t)}{|\omega (x,t)|}
$$
defined in the region $\{x| |\omega(x,t)| >0\}$ satisfies
$$
\int_0^t\int_{\{x| |\omega(x,s)| >0\}}|\omega (x,s)|\left |\nabla \xi (x,s)\right |^2 dx ds \le \frac{1}{2}\nu^{-2}\int_{{\mathbf R}^3}|u_0(x)|^2dx.
$$
\end{thm}
The proof of this result starts with the vorticity equation for the mollified equation:
\begin{equation}
\left (\partial_t + [u]\cdot\nabla - \Delta\right )\omega_i + \epsilon_{ijk}
\partial_j([u]_l)(\partial_l u_k) = 0.\label{omegamol}
\end{equation} 
Here $\epsilon_{ijk}$ is the signature of the permutation $(1,2,3)\mapsto (i,j,k)$ and repeated indices are summed. Multiplying scalarly by $\xi$ one obtains
$$
\left (\partial_t + [u]\cdot\nabla - \Delta\right )|\omega| + \nu|\omega ||\nabla\xi|^2 + Det(\xi, \nabla [u_l], \partial_l u) = 0
$$
where $Det(a,b,c)$ is the determinant of the matrix formed by the three vectors
$a, b, c$.
Integrating  in space and using the energy dissipation one can deduce
the bounds for $\omega$ in $L^1$ and the bound on the 
direction $\xi$. For the bound on $u$ in $L^{\infty}$ one uses the enstrophy
differential inequality
$$
\frac{d}{2dt}\int_{{\mathbf R}^3}|\omega(x,t)|^2 + \nu\int_{{\mathbf R}^3}
|\nabla\omega (x,t)|^2dx \le 
$$
\begin{equation}
\sqrt{4\pi}\left (\int_{{\mathbf R}^3}|\omega (x,t)|^2dx \right)^{\frac{3}{4}}\left (\int_{{\mathbf R}^3}
|\nabla \omega (x,t)|^2dx\right)^{\frac{3}{4}}. \label{ensmol}
\end{equation}
This is obtained from the vorticity equation (\ref{omegamol}) 
above by multiplication by $\omega $, integration by parts, and use of bound
(see (\ref{interpogev}) below)
\begin{equation}
\|u\|_{L^{\infty}} \le \sqrt{4\pi} \|\omega\|_{L^2}^{\frac{1}{2}}\|\nabla\omega\|_{L^2}^{\frac{1}{2}}.\label{interpol}
\end{equation}
Then one employs the idea of \cite{gft}: one divides by
$(c^2 + \int_{{\mathbf R}^3}|\omega (x,t)|^2 dx)^2$, and integrates in time using
the energy principle, ($c>0$ is a constant). One obtains a bound for
$$
\int_0^t\left \{ \int_{{\mathbf R}^3}|\nabla\omega (x,s)|^2dx\right \}^{\frac{1}{3}}ds
$$
in terms of the initial data. The $L^{\infty}$ bound follows from interpolation (\ref{interpol}), the bound above, and the energy principle.  We omit 
further details.

\section{Non-uniform bounds}
If  $u$ is a smooth solution of the Navier-Stokes equations and if  
\begin{equation}
\int _0^T \left (\int_{{\mathbf R}^3} |\omega (x,t)|^2 dx \right )^2dt \le D^3<\infty
\label{suf1}
\end{equation}
then one can bound  any derivative of $u$  on $(0, T]$ in terms of the initial data, viscosity, $T$ and $D$ (see \cite{cfbook}, and references therein). Similarly,  if one has a bound
\begin{equation}
\int_0^T \|u(\cdot, t)\|^2_{L^{\infty}(dx)}dt \le B<\infty \label{suf2}
\end{equation}
then one can bound  any derivative of $u$  on $(0, T]$ in terms of the initial data, viscosity, $T$ and $B$. The quantities $D$ and $B$ have same dimensional count as viscosity (units of length squared per time). If one has a 
regularization that respects the energy dissipation and one has uniform 
bounds for the corresponding quantities then  one can prove global existence 
of smooth solutions. If any of the two conditions is met then the 
solution is real analytic for positive times. Consider the mollified equation (\ref{nsmol}) at $\delta > 0$. Assuming for instance (\ref{suf1}) one obtains bounds for $\sup_{t\le T}\int_{{\mathbf R}^3}|\omega(x,t)|^2dx$ and for $\int_0^T\int_{{\mathbf R}^3}|\nabla \omega(x,t)|^2dxdt$ in terms
of initial data, $D$ and $T$, directly from (\ref{ensmol}). The 
interpolation (\ref{interpol}) then produces a bound for $B$. Vice-versa, if 
one has the assumption (\ref{suf2}) then one does not uses interpolation when one derives the enstrophy inequality from (\ref{omegamol}); rather, one integrates by parts to reveal $u$ and one uses directly the assumption about $\|u\|_{L^{\infty}}$ to deduce a uniform bound for the maximum enstrophy in the time interval $[0,T]$:
\begin{equation}
\sup_{t\le T}\int_{{\mathbf R}^3}|\omega(x,t)|^2dx \le {\mathcal E}<\infty \label{enstr}
\end{equation}
This allows to bound $D$. In either case, the number ${\mathcal E}$ depends 
on the numbers $D$ (respectively $B$) of assumptions (\ref{suf1}) 
(respectively (\ref{suf2})). By increasing ${\mathcal E}$, if necessary, 
we may assume, without loss of generality the
condition
\begin{equation}
\rho = {\mathcal E}^2 T \nu^{-3} > 1. {\label{cond}}
\end{equation}
This condition reflects the fact that we are not pursuing  
decay estimates.  If no assumption is made then ${\mathcal E}$
depends on $\delta>0$. Once the enstrophy is bounded in time, higher derivatives are bounded using the Gevrey-class method of \cite{FT}. 
\begin{thm} Let $\delta > 0$. Consider solutions of (\ref{nsmol})
with initial data $u_0\in L^2$, $\omega_0 \in L^2$. Assume that one of the inequalities
(\ref{suf1}) or (\ref{suf2}) holds on the interval of time
$[0,T]$. Then there exists a constant $c_0 \in (0,1)$ depending only on the
number $\rho = \rho({\mathcal E}, \nu, T)$ of (\ref{cond}), so that
$$
\sup_{t_0\le t\le T}\int_{{\mathbf R}^3} e^{2\lambda |\xi |}|\widehat{\omega }(\xi, t)|^2d\xi
\le 2{\mathcal E} 
$$ 
holds for all $0< t_0\le T$. Here $\lambda = {\sqrt{\nu T}}\min\left \{\frac{t_0}{T}; c_0\right \}$.
If ${\mathcal E}$ is uniform in $\delta $
as $\delta\to 0$, then
the solution of the Navier-Stokes equations with initial data $u_0$ is
real-analytic and obeys the bound  above.
\end{thm}

Using the fact that
$u = \nabla\times\psi $ with divergence-free stream vector $\psi$ 
one sees easily that 
$$
|\widehat{u}(\xi)| \le \frac{1}{|\xi |}|\widehat{\omega} (\xi)| 
$$
holds pointwise. It is elementary then to check that
\begin{equation}
\int_{{\mathbf R}^3}e^{\lambda |\xi|}|\widehat{u}(\xi)| d\xi \le {\sqrt{4\pi}}\left\{\int_{{\mathbf R}^3}e^{2\lambda |\xi|}|\widehat{\omega}(\xi)|^2 d\xi\right \}^{\frac{1}{4}}\left\{\int_{{\mathbf R}^3}e^{2\lambda |\xi|}|\xi |^2|\widehat{\omega}(\xi)|^2 d\xi\right \}^{\frac{1}{4}}
\label{interpogev}
\end{equation}
holds for any positive $\lambda$. Using the Fourier transform of (\ref{nsmol})
and the inequality above one may follow closely the idea of 
\cite{FT}. One considers the quantity
$$
y(t) = \int_{{\mathbf R}^3}e^{2v(t-s)|\xi|}\left |\widehat{\omega }(\xi,t)\right |^2d\xi
$$
with $v=\sqrt{\frac{\nu}{T}}$ and $0\le s\le t\le T$ and derives a differential
inequality of the form
$$
\frac{dy}{dt} \le c\left (\frac{y}{\nu}\right )^3 + c\frac{v^2}{\nu}y
$$
with absolute constant $c$. At $t=s$ one has by construction
$y(s)\le {\mathcal E}$. 
The differential inequality guarantees that $y$ does not exceed $2{\mathcal E}$ on a time interval $s\le t\le s+ 2c_0 T$ with $c_0$ a small 
non-dimensional constant (proportional to  $\rho^{-1}= \nu^3T^{-1}{\mathcal E}^{-2}$). 
The time step $2c_0T$ is uniform because of the assumption that ${\mathcal E}$ is finite. One starts from $s=0$. The differential inequality implies 
a non-trivial Gevrey-class bound for the second half of the first time interval
$[0, 2c_0T]$. One now sets $s=c_0T$ and makes another step of duration $2c_0T$. At each step the second half of the time interval yields a nontrivial bound, and because the initial point $s$ is advanced by a half step one covers all of 
$[c_0 T, T]$. The pre-factor in the exponential bound is uniformly bounded below by $2c_0vT$. If $t_0\le c_0T$ we can advance $2t_0$ at a time, and obtain 
similarly a uniform bound on $[t_0, T]$, with exponent pre-factor bounded below
by $2t_0v$. This result implies that the velocity is real analytic. One may extend $u$ to a complex domain $z =x+iy$ by setting
$$
u(x+iy,t) = (2\pi)^{-3}\int_{{\mathbf R}^3}e^{i(x+iy)\cdot\xi}\widehat{u}(\xi,t)d\xi
$$
and, in view of (\ref{interpogev}), the integral is absolutely convergent for
$|y| \le \lambda $, $t\ge t_0$. Note also that, at fixed $\delta>0$ one has a 
finite, $\delta$-dependent bound on ${\mathcal E}$, and therefore the 
solutions of the mollified equations are real-analytic. Likewise, if one allows for $\delta$ dependent bounds, then the vortex method also can be shown to have real-analytic solutions. In this case, however, because the energy principle
is not strong enough one starts from the assumption (\ref{enstr}):
\begin{thm} Let the initial vorticity $\omega_0$ belong to $L^2$. Consider, for
$\delta >0$ the solution of (\ref{bob}) and assume that (\ref{enstr}) holds
on a time interval $[0, T]$.  Then there exists a constant $c_0 \in (0,1)$ depending only on  the
number $\rho = \rho({\mathcal E}, \nu, T)$ of (\ref{cond}), so that
$$
\sup_{t_0\le t\le T}\int_{{\mathbf R}^3} e^{2\lambda |\xi |}|\widehat{\omega }(\xi, t)|^2d\xi
\le 2{\mathcal E} 
$$ 
holds for all $0< t_0\le T$. Here $\lambda = {\sqrt{\nu T}}\min\left \{\frac{t_0}{T}; c_0\right \}$.
If ${\mathcal E}$ is uniform in $\delta $
as $\delta\to 0$, then the solution of the Navier-Stokes equations (\ref{omegu}) with initial data $\omega_0$ is
real-analytic and obeys the bound  above.
\end{thm}
 
The proof follows the same ideas as the proof of the corresponding result for the mollified equation (\ref{nsmol}).

\section{Euler equations}
The three dimensional Euler equations
\begin{equation}
\partial_tu + u\cdot\nabla u + \nabla p = 0,\,\,\,\,\,\,\, \nabla\cdot u = 0 \label{euler}
\end{equation}
are locally well-posed \cite{ema}, \cite{ka}. They 
conserve kinetic energy (if the solutions are smooth enough \cite{o}, \cite {ey}, \cite{cet}). Such smooth solutions can be interpreted 
\cite{a} as geodesic paths on an infinite dimensional group of 
transformations. Despite of energy conservation, gradients
of solutions may grow \cite{serre}. The 
vorticity $\omega = \nabla\times u$  obeys
\begin{equation}
\left (\partial _t + u\cdot\nabla \right )\omega = \omega\cdot\nabla u.\label{omege}
\end{equation}
Because of the quadratic nature of this equation and the fact that
the strain matrix 
\begin{equation}
S = \frac{1}{2}\left \{\left (\nabla u)\right ) + \left (\nabla u\right )^*
\right \}\label{strain}
\end{equation}
is related to the 
vorticity by a linear classical singular Calderon-Zygmund  
integral, it was suggested \cite{Ma} that blow up of the vorticity might 
occur in finite time. This problem is open, despite much research \cite{child}, \cite{cnote}, \cite{cimr}, \cite{clm}, \cite{go}, \cite{jt}.
The blow up cannot occur unless the time integral of the maximum modulus
of vorticity diverges \cite{bkm}. The vorticity magnitude obeys
\begin{equation}
\left (\partial_t + u\cdot\nabla\right )|\omega | = \alpha|\omega|,\label{stretch}
\end{equation}
where the logarithmic material stretching rate $\alpha$ can be represented 
\cite{siam} as
\begin{equation}
\alpha (x, t) = \frac{3}{4\pi} P.V.\int D(\widehat{y}, \xi(x-y, t), \xi(x,t))
|\omega(x-y,t)|\frac{dy}{|y|^3}.\label{alpha}
\end{equation} 
Here $\widehat{y} = \frac{y}{|y|}$,  
\begin{equation}
D(\widehat{y}, \xi(x-y, t), \xi(x,t)) = \left (\widehat{y}\cdot\xi(x,t)\right )Det( \widehat{y}, \xi(x-y, t), \xi(x,t))\label{geomfact}
\end{equation}
and
\begin{equation}
\xi (x,t) = \frac{\omega(x,t)}{|\omega (x,t)|}.\label{vordir}  
\end{equation}
In the two dimensional case $\xi = (0,0,1)$ and so $\alpha = 0$. In three 
dimensions, if the vorticity direction is well-behaved locally in regions of 
high vorticity, then there is a geometric depletion of nonlinearity. More precisely,
if $\xi(x,t)\times \xi(x-y,t)$ vanishes in a 
quantitatively controlled fashion as $y\to 0$, (for instance $|\xi (x-y,t)\times \xi (y,t)| \le k|y|$),  then $\alpha $ can be bounded 
in terms of less singular integrals (for instance in terms of velocity instead of vorticity). This observation suggests a correlation between vorticity growth and the geometry of vortex tubes. Such a correlation has been observed in numerical studies and 
was exploited to prove conditional results regarding blow up for the
Euler equations, \cite{siam}, \cite{cfm},  the Navier-Stokes equations \cite{cfiu}, 
and the quasi-geostrophic model \cite{Co}. The quasi-geostrophic model \cite{siam},
 \cite{ccw}, \cite{CMT}, \cite{CW1}, \cite{R} is an example of an active scalar.
Active scalars are advection-diffusion evolution equations for scalar quantities advected by an incompressible velocity they create: the velocity is obtained from the scalar using a fixed, time-independent formula:
\begin{equation}
\left (\partial_t + u\cdot\nabla -\kappa\Delta \right )q = 0\,\,\,\, u = U[q].\label{model}
\end{equation}
The Euler
equations themselves are an active vector system \cite{cel}, \cite {celw}:
\begin{equation}
\left (\partial_t + u\cdot\nabla \right )A = 0, \label{Au} 
\end{equation}
with
\begin{equation}
u = W[A]\label{Wau}
\end{equation}
the Weber formula \cite{serr} 
\begin{equation}
u(x,t) = {\mathbf P}\left ((\nabla A)^*v\right).\label{web}
\end{equation} 
Here ${\mathbf P}$ is the Leray-Hodge projector on divergence-free
functions, and 
\begin{equation}
v(x,t) = u_0(A(x,t))\label{va}
\end{equation}
is a solution of
\begin{equation}
\left (\partial_t + u\cdot\nabla \right )v = 0\label{vpass}
\end{equation}
The initial data for $A$ is the identity
\begin{equation}
A(x,0) = x, \label{idem}
\end{equation}
and the initial datum for $v$ is $u_0$. Thus $A(x,t) =a$ is the inverse of the Lagrangian
path $a\mapsto X(a,t)$. The familiar Cauchy formula is in this language 
\begin{equation}
 \omega_q = \frac{1}{2}\epsilon_{qij}\left (Det\left
[\zeta ; \frac{\partial A}{\partial x_i}; \frac{\partial A}{\partial
  x_j}\right ]\right ). \label{cauch}
\end{equation}
Here 
\begin{equation}
\zeta (x,t) = \omega_0(A(x,t))\label{zeul}
\end{equation}
is the solution of
\begin{equation}
\left (\partial_t + u\cdot\nabla \right )\zeta = 0\label{znul}
\end{equation}
with initial datum $\zeta (x,0) = \omega_0(x)$.   

One may use a near-identity approach to the incompressible Euler equations:
One solves (\ref{Au}) for short time, as long as $\nabla A - I$ is not too large.
Then one stops and resets $u =W[A]$ in place of $u_0$, sets  $A =x$, adjusts the clock, and starts again.  This approach allows one to interpret the condition \cite{bkm}
of absence of blow up in terms of  $\nabla A$ : if
$$
\int_0^T\|\nabla A(\cdot , t)\|^2_{L^{\infty}}dt <\infty
$$
then the solution of the Euler equations is smooth \cite{cel}.   

\section{Diffusive Lagrangian transformations}
The central object in  the Lagrangian description of fluids  is
the Lagrangian path transformation $a\mapsto X(a,t)$; $x =X(a,t)$ represents 
the position at time $t$ of the fluid particle that started at $t=0$ from
$a$. At time $t=0$ the transformation is the identity, $X(a,0) = a$.  
An Eulerian-Lagrangian formulation of the Navier-Stokes equations 
\cite{xlns} parallels the active vector formulation of the Euler equations \cite{cel}.  In order to unify the exposition we associate to a 
given divergence-free velocity $u(x,t)$ the 
operator
\begin{equation}
\partial_t + u \cdot\nabla -\nu\Delta = {\Gamma_{\nu} (u,\nabla )}\label{L}.
\end{equation} 
Using $\Gamma_{\nu}(u,\nabla)$ we associate to any 
divergence-free, time dependent velocity field
$u$ a transformation $x\mapsto A(x,t)$ that obeys
\begin{equation}
\Gamma _{\nu}(u, \nabla) A = 0, \label{aeq}
\end{equation}
with initial data
\begin{equation}
A(x,0) = x\label{azero}
\end{equation}
Boundary conditions are imposed by considering the displacement
vector 
\begin{equation}
\ell(x,t) = A(x,t) - x\label{d}
\end{equation}
that joins the Eulerian position $x$ to the the diffusive label $A$. This is
required to vanish at infinity, and one can think of $A$ as being computed by
solving  
\begin{equation}
\left (\partial_t   + u\cdot\nabla  -\nu \Delta \right )\ell 
+ u = 0
\label{deltaeq}
\end{equation}
with initial data 
\begin{equation}
\ell (x,0) = 0.\label{deltain}
\end{equation}
If $u$ is the solution of the Euler equation and $\nu = 0$ in the equations
(\ref{aeq}, \ref{deltaeq}) then the map $A$ is the inverse of the 
particle trajectory map $a\mapsto x= X(a,t)$. In the presence of viscosity
this map obeys a diffusive equation, departing thus from its conventional
interpretation as inverse of particle trajectories. Nevertheless, continuing 
the analogy with the inviscid situation, one uses the map 
$x\mapsto A(x,t)$ to pull back the Lagrangian differentiation with respect 
to particle position and write it in Eulerian coordinates. This 
Eulerian-Lagrangian derivative is given by
\begin{equation}
\nabla^A = Q^*\nabla \label{la}
\end{equation}
where
\begin{equation}
Q(x,t) = \left (\nabla A(x,t)\right )^{-1},\label{Q}
\end{equation}
that is
\begin{equation}
\nabla^A_i = Q_{ji}\partial_j.\label{lai}
\end{equation} 
In the case $\nu = 0$ the invertibility of $\nabla A$ follows from incompressibility; in the diffusive case the determinant of $\nabla A$ does not remain identically equal to one as time passes. This imposes a constraint on the time 
of integration.  We consider a small non-dimensional parameter $g>0$ and
work with the constraint
\begin{equation}
\sup_{0\le t\le \tau}\sup_{x\in {\mathbf R}^3}|\nabla \ell (x,t)| \le g.
\label{g}
\end{equation}
With this constraint satisfied we can guarantee the invertibility of 
$\nabla A$. In order to describe the dynamics and their relationship to the 
Eulerian dynamics one needs to consider second derivatives of $A$. These influence
the dynamics because commutators between Eulerian-Lagrangian and 
Eulerian derivatives do not vanish, in general:
\begin{equation}
\left [\nabla^A_i, \nabla_k\right ] = C^{m}_{k; i}\nabla^A_m.\label{str}
\end{equation}
The coefficients $C^{m}_{k;i}$ are given by 
\begin{equation}
C^{m}_{k;i} = \{\nabla^A_{i}(\partial_k \ell_m)\}.\label{cmki}
\end{equation}
Note that 
$$
C^{m}_{k;i} = Q_{ji}\partial_j\partial_k A^m =
\nabla^A_i(\partial_k A^m) = [\nabla^A_i,\partial_k]A^m.
$$
These commutator coefficients are related to the Christoffel coefficients of 
the usual flat connection in ${\mathbf R}^3$ computed using the change of variables $a= A(x,t)$. With this change of variables, a  straight line 
in $x$, $x(s) = ms + b$ becomes the label path $a(s) = A(x(s),t)$ and the 
geodesic equation $\frac{d^2 x}{ds^2} = 0$ becomes
$\frac{d^2 a^m}{ds^2} + \Gamma^m_{ij}\frac{da^i}{ds}\frac{da^j}{ds} = 0$ with
$$
\Gamma^m_{ij} = -C^{m}_{k;j}Q_{ki}.
$$ 
The simple geometry of ${\mathbf R}^3$ is hidden behind a complicated 
transformation, but the transformation is the main object of study.  
The coefficients $C^{m}_{k;i}$ (but not $u$) enter the 
commutation relation between
the Eulerian-Lagrangian label derivative and $\Gamma_{\nu}(u, \nabla)$:
\begin{equation}
\left [{\Gamma_{\nu} (u, \nabla)}, \nabla^A_{i}\right ] = 2\nu C^{m}_{k ;i}
\partial_k \nabla^A_{m}
 \label{lagcom}
\end{equation}
This commutation relation is the viscous counterpart of the inviscid 
commutation of time and label derivatives.

In the inviscid case the map $A$ is the main active ingredient in the dynamics.
The Weber formula (\ref{web}) \cite{serr} computes the velocity at time 
$t$ directly from the gradient of $A$ using a passively advected velocity 
$v$ (\ref{va}). In the viscous case, the Weber formula
\begin{equation}
u = {\mathbf P}\left ((\nabla A)^*v\right )\label{pu}
\end{equation}
can still be used but $v(x,t)$ is no longer passive. Instead of (\ref{vpass}) 
$v$ obeys
\begin{equation}
\Gamma_{\nu}(u,\nabla ) v = 2\nu C\nabla v, \label{psieq}
\end{equation}
that is
\begin{equation}
\Gamma_{\nu}(u,\nabla)v_i = 2\nu C^{m}_{k;i}\partial_kv_m \label{veq}
\end{equation}
with initial data 
\begin{equation}
v(x,0) = u_0(x).
\label{psin}
\end{equation}
The equations (\ref{aeq}), (\ref{psieq}) together with the Weber formula
(\ref{pu}) are equivalent to the Navier-Stokes equations \cite{xlns}:
\begin{thm}
Let  $A$, $v$ and $u$ solve the system
(\ref{aeq},\,\ref{psieq},\,\ref{pu}). Then $u$ obeys the 
incompressible Navier-Stokes equations, 
$$
\partial_t u + u\cdot\nabla u - \nu\Delta u + \nabla p = 0, \quad \nabla\cdot u =0.
$$ 
\end{thm}
We describe now kinematic consequences of (\ref{aeq}). One starts with
a smooth arbitrary incompressible velocity field
$u$, one computes $A$ using (\ref{aeq}), then one computes the inverse 
matrix $Q= (\nabla A)^{-1}$ (at least for a short time), and then 
Eulerian-Lagrangian derivatives 
$\nabla^A$ and  coefficients $C$. One then can evolve a vector $v$ from an initial datum solving (\ref{psieq}), and
one can compute its Eulerian-Lagrangian curl $\zeta = \nabla^A\times v$. 
The resulting equations \cite{xlns}, \cite{cfil} are summarized below.
 
\begin{thm}
Let $u$ be an arbitrary divergence-free function and associate to it
a map $A$ solving (\ref{aeq}) and a vector field $v$ solving (\ref{psieq}).
Then $w$ defined by 
\begin{equation}
w_i = (\partial_iA^m)v_m
\label{w}
\end{equation}
obeys the cotangent equation
\begin{equation}
{\Gamma_{\nu}(u,\nabla) }w + (\nabla u)^*w = 0.\label{weq}
\end{equation}
The Eulerian curl of $w$, $\nabla\times w$ obeys the equation
\begin{equation}
\Gamma_{\nu}(u, \nabla )\left (\nabla\times w\right ) = (\nabla u)(\nabla\times w).
\label{nablaweq}
\end{equation}
The Eulerian-Lagrangian curl of $v$, $\zeta = \nabla^A\times v$, obeys
\begin{equation}
\Gamma_{\nu}(u,\nabla) \zeta_q = 2\nu C^{m}_{k;m}\partial_k\zeta_q - 
2\nu C^{q}_{k;j}\partial_k\zeta_j + \nu C^{m}_{k;i}C^{r}_{k;j}
\epsilon_{qji}\epsilon_{rmp}\zeta_p.\label{zetaueq}
\end{equation}
The Eulerian curl of w, $\nabla \times w$ and Eulerian-Lagrangian 
curl of $v$, $\zeta = \nabla^A\times v$ are related by the formula
\begin{equation}
\left (\nabla_E\times w \right )_q =
\frac{1}{2}\epsilon_{qij}\left (Det\left
[\zeta ; \frac{\partial A}{\partial x_i}; \frac{\partial A}{\partial  x_j}\right ]\right ) 
\label{omegauzeta}
\end{equation}
The determinant of $\nabla A$ obeys
\begin{equation}
\Gamma_{\nu}(u,\nabla )\left (\log Det(\nabla A)\right )
= \nu\left\{C^{i}_{k;s}C^{s}_{k;i}\right\}.\label{deteq}
\end{equation}
\end{thm}

These considerations apply to arbitrary
$u$ without having to impose the equation of state (\ref{pu}).
When $\nu = 0$ the cotangent equation (\ref{weq}) appears in 
Hamiltonian formalisms \cite{kuz}, \cite{ose} for the Euler equation
in various gauges. The numerical merits of these  have been analyzed 
critically \cite{russo}. When $u$ is related to $w$ by a
filtered Weber formula, a gauge of the cotangent equation appears as 
a model of formally ``averaged'' Euler equation \cite{mar} and, 
in a viscous case, as a model of Reynolds' equation \cite{holm}.

Note that, when $\nu = 0$, the equation (\ref{zetaueq}) is just
the pure advection equation 
$\left (\partial_t + u\cdot\nabla \right )\zeta = 0$. If $\nu>0$, 
$\zeta$ obeys a linear dissipative equation with coefficients $C^{m}_{k;i}$. 
Using Schwartz inequalities only, one obtains
\begin{equation}
\Gamma_{\nu}(u, \nabla) |\zeta(x,t) |^2 + \nu |\nabla \zeta (x,t) |^2 
\le 17 \nu |C(x,t)|^2|\zeta (x,t)|^2
\label{zetaineq}
\end{equation}
where 
$$
|C(x,t)|^2 = C^{m}_{k;i}(x,t)C^{m}_{k;i}(x,t),\quad |\zeta(x,t)|^2 = \zeta_q(x,t)\zeta_q(x,t)
$$
are squares of Euclidean norms.

The evolution of the coefficients $C^{m}_{k;i}$
defined in (\ref{cmki}) is given by
$$
{\Gamma_{\nu}(u,\nabla)}\left (C^{m}_{k;i} \right ) =
 - (\partial_lA^m)\nabla^A_i(\partial_k(u_l))
$$
\begin{equation}
-(\partial_k(u_l))C^{m}_{l;i}+ 2\nu C^{j}_{l;i}\cdot\partial_l\left 
(C^{m}_{k;j}\right ).\label{curveq}
\end{equation}
At time $t=0$ the coefficients vanish, $C^m_{k;i}(x, 0) = 0$.
Note that the linear equation (\ref{nablaweq}) is 
identical to the nonlinear vorticity equation (\ref{omegu}).
It is the relation (\ref{pu}) that decides whether or not we are
solving the Navier-Stokes equation; if $u$ is a solution of the
Navier-Stokes equations then (\ref{pu}) means $u = {\mathbf P}w$ and
consequently $\nabla\times u = \nabla \times w$ implies that the 
Eulerian curl of $w$ is the fluid's vorticity. The formula
(\ref{omegauzeta}) relating the vorticity to $\zeta$ is then a 
viscous counterpart of the  Cauchy formula.  

\begin{thm} If $u$ solves the Navier-Stokes equation,
$A$ solves (\ref{aeq}), and $v$ solves (\ref{psieq}) then 
the Eulerian curl of $u$, $\omega = \nabla\times u$ is related
to the Eulerian-Lagrangian curl of $v$, $\zeta = \nabla^A\times v$
by the Cauchy formula
\begin{equation}
\omega_q = \frac{1}{2}\epsilon_{qij}\left (Det\left
[\zeta ; \frac{\partial A}{\partial x_i}; \frac{\partial A}{\partial
  x_j}\right ]\right ). \label{omegazeta}
\end{equation}
\end{thm} 

Because of the linear algebra identity
$$
((\nabla A)^{-1}\zeta)_q  = \left (Det(\nabla A)\right )^{-1}\frac{\epsilon_{qij}}{2}\left (Det{\left[\zeta ; \frac{\partial A}{\partial x_i}; \frac{\partial A}{\partial x_j}\right ]}\right )
$$
one has
\begin{equation}
\omega = \left (Det(\nabla A)\right )(\nabla A)^{-1}\zeta.\label{cauchynu}
\end{equation}
In  two-dimensions (\ref{omegazeta}, \ref{cauchynu}) become
\begin{equation}
\omega = \left (Det(\nabla A)\right )\zeta,\label{twodcauchynu}
\end{equation}
reflecting the fact that, for $\nu = 0$, $\omega = \zeta$ in that case.
A consequence of (\ref{omegazeta}) or ({\ref{cauchynu}) is the identity
\begin{equation}
\omega\cdot\nabla_E = (Det(\nabla A))\left (\zeta\cdot\nabla_A\right )\label{derivs}
\end{equation}
that generalizes the corresponding inviscid identity. Let us consider the 
expression
\begin{equation}
{\mathcal  C}(q,M) = (Det M)M^{-1}q \label{Ca}
\end{equation}
defined for any pair $(q, M)$, where $q\in {\mathbf R}^3$, 
$M\in GL({\mathbf R}^3)$ are, respectively, a vector and an 
invertible matrix. This expression, underlying the Cauchy formula, is 
linear in $q$ and quadratic in $M$,
\begin{equation}
{\mathcal C}(q,M) _k = \frac{1}{2}\epsilon_{ijk}Det \left (M_{.\,,i}, M_{.\,, j}, q\right ).
\label{qua}
\end{equation}
The quadratic expression in the right hand side is 
defined for any matrix $M$. It is easy to check that
\begin{equation}
{\mathcal C}(q,MN) = 
{\mathcal C}\left ({\mathcal C}(q,M) , N\right) \label{cyc}
\end{equation}
and 
\begin{equation}
{\mathcal C}(q, {\mathbf  I})=  q\label{idd}
\end{equation}
hold, so ${\mathcal C}$ describes an action of $GL({\mathbf R}^3)$ in  $R^3$.
A third property follows from the explicit quadratic expression (\ref{qua})
\begin{equation}
{\mathcal C}(q, {\mathbf 1} + N)  =  (1+ Tr(N))q - Nq + 
{\mathcal C}(q,N)\label{nid}
\end{equation}
Here $N$ is any matrix, and the meaning of ${\mathcal C}(q,N)$ is given 
by (\ref{qua}).
If we consider, instead of vectors $q$ and matrices $M$, vector valued 
functions $q(x,t)$ and matrix valued $M(x,t)$ and use the same formula
\begin{equation}
{\mathcal C}(q, M)(x,t)  =  {\mathcal C}(q(x,t), M(x,t))  = (Det(M(x,t))
\left (M(x,t)\right )^{-1}q(x,t)
\label{bigC}
\end{equation}
then the  properties  (\ref{cyc}, \ref{idd}, \ref{nid}) as well 
as (\ref{qua}) obviously still  hold.

Denote
$$
\epsilon (s) = \nu \int |\nabla u (x,s)|^2dx
$$
the total instantaneous energy dissipation rate,
and
$$
K(t) =  \frac{1}{2}\int|u(x,t)|^2dx
$$
the  total kinetic energy. We consider functions $u(x,t)$ that satisfy
\begin{equation}
K(t) + \int_0^t \epsilon (s) ds \le
K_0,  
\label{epsbound}
\end{equation}
and
\begin{equation}
\int_{0}^t\|u(\cdot,s)\|_{L^{\infty}(dx)}ds \le K_{\infty}\label{maxu}
\end{equation}
For solutions 
of the Navier-Stokes equations and for solutions of the mollified
equation (\ref{nsmol})  the constants $K_0$, $K_{\infty}$ depend
on the initial kinetic energy, viscosity and time only. For vortex methods, however,
these constants depend on the cut-off scale $\delta $. 
The displacement $\ell$ satisfies certain bounds that follow from
the bounds above and (\ref{deltaeq}). 
\begin{thm} Assume that the vector valued function $\ell$ obeys
(\ref{deltaeq}) for $t\in [t_1, T]$ and $\ell(\cdot, t_1) = 0$, for some $t_1\ge 0$. Assume that the velocity $u(x,t)$ 
is a divergence-free  function that satisfies 
the bounds (\ref{epsbound}) and (\ref{maxu}) on the time interval $[0, T]$. Then $\ell$ satisfies
the inequality 
\begin{equation}
\|\ell (\cdot, t) \|_{L^{\infty}} \le K_{\infty}\label{maxdel}
\end{equation}
together with
\begin{equation}
\int|\ell (x,t)|^2 dx \le  
2K_0(t-t_1)^2, 
\label{ltwo}
\end{equation}
\begin{equation}
\int_{t_1}^t\int |\nabla \ell(x,s)|^2 dx ds \le 
\frac{K_0(t-t_1)^2}{\nu},
\label{nablaeltwo}
\end{equation}
and
\begin{equation}
\int |\nabla \ell(x,t)|^2 dx + \nu \int_{t_1}^t\int |\Delta\ell(x,s)|^2dxds\
\le C_1\left (\frac{K_0 (t-t_1)}{\nu} +
\frac{K_{\infty}^2E_0}{\nu^2}\right ).
\label{deltaltwo}
\end{equation}

\end{thm}

Let us consider the analytic norms 
\begin{equation}
\|u\|_{\{A,\lambda, p\}} = \left \{\int e^{p \lambda |\xi|} |\widehat{u}(\xi)|^pd\xi\right\}^{\frac{1}{p}}\label{anorms}
\end{equation}
We will use $p =1$ and $p=2$. One can prove 

\begin{thm}
Let the vector valued function $\ell$ solve
(\ref{deltaeq}) for $t\ge t_1\ge 0$ with $\ell(\cdot ,t_1) = 0$.
Assume that a velocity $u(x,t)$ defined for $t\ge t_0$ , $t_0\le t_1$ 
is a divergence-free  function that satisfies
\begin{equation}
\sup_{t_0\le t\le T} \|u\|_{\{A,r,1\}} \le U_r.\label{bou}
\end{equation}
Then there exists an absolute constant (a pure number) $c$ such that
\begin{equation}
\|\ell (\cdot, t )\|_{\{A, r,1\}} \le (t-t_1)U_re^{c\frac{(t-t_1)U_r^2}{\nu}}\label{lbound}
\end{equation}
holds for $t\ge t_1$.
\end{thm}

Note that (\ref{interpogev}) reads
\begin{equation}
\|u\|_{\{A, \lambda, 1\}} \le \sqrt{4\pi}\|\omega\|_{\{A, \lambda, 2\}}^{\frac{1}{2}}\|\nabla\omega\|_{\{A, \lambda, 2\}}^{\frac{1}{2}}
\label{interlambedn}
\end{equation}
Note also that, if $r<\lambda$ then
\begin{equation}
\|\nabla\omega\|_{\{A, r, 2\}} \le \frac{1}{e(\lambda -r)}\|\omega\|_{\{A, \lambda, 2\}}\label{updown}
\end{equation}
and
\begin{equation}
\|\nabla\omega\|_{\{A, r, 1\}} \le \frac{1}{e(\lambda -r)}\|\omega\|_{\{A, \lambda, 1\}}\label{updown1}
\end{equation}
hold. Combining Theorem 3 or Theorem 4 with the preceding result we obtain
therefore

\begin{thm} 
Consider solutions of the mollified equations (\ref{nsmol}) 
or of the vortex method (\ref{bob}) associated to a filter (\ref{[u]}).
Assume the initial  data are divergence-free and belong to $H^1$, 
$u_0\in L^2$, $\nabla\times u_0 =\omega_0\in L^2$. Consider $0\le t\le T$
and let ${\mathcal E}$ denote a bound for the enstrophy on the time interval
$[0,T]$
$$
\sup _{0\le t\le T}\int_{{\mathbf R}^3}|\omega (x,t)|^2dx \le {\mathcal E}.
$$ 
Consider an arbitrary transience time $0<t_0<T$ and the length
scale
$$
\lambda = \sqrt{\frac{\nu}{T}}\min\left \{t_0;\, c_1\nu^3{\mathcal E}^{-2} \right \}
$$
with $c_1$ a certain absolute constant. Then the velocity $u$ obeys the bound 
(\ref{bou}) for any $r<\lambda$, with
$U_r$ given by
\begin{equation}
U_r = c_2(\lambda -r)^{-\frac{1}{2}}{\mathcal E}^{\frac{1}{2}}\label{ur}
\end{equation}
and $c_2 = 2\sqrt{\frac{2\pi}{e}}$. Consequently, for any $t_1\ge t_0$, the 
solution $\ell$ of the equation (\ref{deltaeq}) with initial data 
$\ell (\cdot, t_1) = 0$ obeys the bound (\ref{lbound}) and, for 
arbitrary $r_1<r<\lambda$,
\begin{equation}
\|\nabla \ell(\cdot, t)\|_{\{A, r_1, 1\}}\le e^{-1}(r-r_1)^{-1}(t-t_1)U_re^{c\frac{(t-t_1)U_r^2}{\nu}}\label{grl}
\end{equation}
holds for $t\ge t_1$.
\end{thm}

These bounds depend on the cut-off scale $\delta$ of the 
filter only through the bound on enstrophy ${\mathcal E}$. In particular, 
if the enstrophy is bounded uniformly for small $\delta$ 
on a time interval then the results above apply for the 
Navier-Stokes equations on that time interval.  
If we measure length in units of $\sqrt{\nu T}$ and time in units of
$T$ then the enstrophy bound for $u(x,t) = \sqrt{\frac{\nu}{T}}\tilde{u}\left (\frac{x}{\sqrt{\nu T}}, \frac{t}{T}\right )$ becomes
$$
\sup_{0\le s\le 1} \|\nabla \tilde u(\cdot,  s)\|_{L^2} \le G
$$
with $G$ the non-dimensional number given by 
\begin{equation}
G^2 = \nu^{-\frac{3}{2}}T^{\frac{1}{2}}{\mathcal E}\label{g2}
\end{equation}
Note that the number $\rho$ of (\ref{cond}) is just $\rho = G^4$. In these
units $G$ is the only solution dependent parameter that we do not control. 
In terms of $G$ and in these units, the definition of $\lambda $ becomes
\begin{equation}
\frac{\lambda}{\sqrt{\nu T}} = \tilde{\lambda} = \min\{s_0; \,c_1G^{-4}\},\label{latilde}
\end{equation}
with $s_0\in (0,1)$ arbitrary. The bound (\ref{ur}), with
for $0<\tilde{r}<\tilde{\lambda}$, $r =\sqrt{\nu T}\tilde{r}$, becomes
\begin{equation}
\sqrt{\frac{T}{\nu}}U_r = \tilde{U}_{\tilde r} \le c_2 (\tilde{\lambda}-\tilde{r})^{-\frac{1}{2}} G.\label{utilde}
\end{equation}
For fixed $s_0$ and large $G$ we may take $\tilde{\lambda} \sim G^{-4}$
so that if we consider $\tilde{r} = (1-\gamma)\tilde{\lambda}$ with
$0<\gamma<1/4$ we deduce 
\begin{equation}
\tilde {U}_{\tilde{r}} \le c_3G^3 \label{tildu}
\end{equation}
Choosing $r_1 = (1-\gamma)r$ in (\ref{grl}) we deduce from (\ref{grl})
\begin{equation}
\|\nabla \ell (\cdot, t)\|_{\{A; r_1, 1\}} \le g\label{ggrl}
\end{equation}
for $t_0\le t_1\le t\le t_1 +\tau T$ with
\begin{equation}
\tau = c_4g G^{-7}.\label{tau}
\end{equation}
Choosing $r_2 = (1-\gamma)r_1$ we deduce that
\begin{equation}
\|\nabla\nabla \ell (\cdot, t)\|_{\{A; r_2, 1\}} \le c_5(\nu T)^{-\frac{1}{2}}G^4 g
\label{nagrl}
\end{equation}
holds on the same time interval.

\section{Conclusions} The viscous Navier-Stokes equations and their
approximations can be described using diffusive, near-identity 
transformations. The velocity is obtained from the near-identity
transformation using a Weber formula and a virtual velocity. The vorticity
is obtained from the near-identity transformation using a Cauchy formula
and a virtual vorticity. The virtual velocity and the virtual vorticity
obey diffusive equations, which reduce to passive advection formally,
if the viscosity is zero. Apart from being proportional to the viscosity,
the coefficients of these diffusion equations involve second derivatives
of the near identity transformation and are related to the Christoffel 
coefficients. If and when the near-identity transformation departs
excessively from the identity, one resets the calculation. Lower bounds
on the minimum time between two successive resettings are given in 
terms of the maximum enstrophy.

\noindent{\bf Acknowledgment.} Research partially supported by NSF-DMS9802611.


\begin{thebibliography}{99}
\bibitem{a} V. I. Arnol'd, {\em Mathematical methods of classical
mechanics}, Springer-Verlag, GTM {\bf 60} New York, 1984.

\bibitem{bkm} J. T. Beale, T. Kato, A. Majda,  Remarks on the breakdown
of smooth solutions for the 3-D Euler equations, {\em Commun. Math. Phys.} {\bf 94} (1984), 61-66
\bibitem{ckn} L. Caffarelli, R. Kohn, L. Nirenberg, Partial regularity for suita ble weak solutions of the Navier-Stokes equations, {\em Commun. Pure Appl. Math.},
{\bf 35} (1982), 771-831.
\bibitem{holm} Shiyi Chen, Ciprian Foias, Darryl D Holm, Eric Olson, 
Edriss S Titi, Shannon Wynn, A connection between the
Camassa-Holm equations and 
turbulent flows in channels and pipes, {\em Phys. Fluids}, {\bf{11}}
(1998), 2343-2353

\bibitem{child} S. Childress, G. R. Ierley, E. A. Spiegel, and W.R. Young, 
  Blow-up of unsteady two-dimensional Euler and Navier-Stokes solutions having
  stagnation-point form, {\em J. Fluid Mechanics} {\bf 203}, 1989, 1 - 22.
\bibitem{chorin} A. Chorin, Numerical study of slightly viscous flow, {\em 
J. Fluid. Mech.} {\bf 57} (1973), 785 - 796.


\bibitem{cnote} P. Constantin, Note on loss of regularity for solutions of the 3D incompressible Euler and related equations, {\em Commun. Math. Phys.} {\bf 106}
(1986), 311 - 326.
\bibitem{c91} P. Constantin,  Navier-Stokes equations and area of interfaces, {\em Commun. Math. Phys.} {\bf 129} (1990), 241 - 266. 
\bibitem{siam} P. Constantin,  Geometric statistics in turbulence, {\em Siam Review}, {\bf 36}, (1994), 73-98. 
\bibitem{cnot} P. Constantin, A few results and open problems regarding incompressible fluids, {\em Notices of the AMS} {\bf 42} (1995), 658-663.
\bibitem{cimr} P. Constantin, The Euler Equations and Nonlocal Conservative Riccati Equations, {\em Intern. Math. Res. Notes}, {\bf 9} (2000), 455-465 .
\bibitem{cel} P. Constantin, An Eulerian-Lagrangian approach for incompressible fluids: Local theory,  {\em J. Amer. Math. Soc.} {\bf 14} (2001), 263-278..
\bibitem{celw} P. Constantin, An Eulerian-Lagrangian approach to fluids, e-print (1999) http://www.aimath.org. 
\bibitem{xlns} P. Constantin, An Eulerian-Lagrangian approach for the
Navier-Stokes equations,  
{\em Commun. Math. Phys.} {\bf{216}} (2001), 663-686.
\bibitem{cfil} P. Constantin, Filtered Viscous Fluid Equations, {\em 
Computer and Mathematics with Applications}, to appear (2001)
\bibitem{c2k} P. Constantin, Some open problems and research directions 
in the mathematical study of fluid dynamics, in {\it Mathematics unlimited:
2001 and beyond}, B. Engquist and W. Schmid edtrs, Springer, Berlin-New York,
2001.
\bibitem{ccw} P. Constantin, D. Cordoba, J. Wu, On the critical
dissipative quasi-geostrophic equation, {\em Indiana Univ. Math. J.}{\bf 50}
(2001), 97-107.
\bibitem{cet} P. Constantin, W.  E, E. Titi, 
Onsager's conjecture on the energy conservation for
solutions of Euler's equation, {\em Commun. Math.
Phys.}, {\bf 165} (1994), 207-209.
\bibitem{cfiu} P. Constantin, C. Fefferman, Direction of vorticity and the problem of global regularity for the Navier-Stokes equations, {\em 
Indiana Univ. Math. J.}, {\bf 42} (1993), 775-787.
\bibitem{cfm} P. Constantin, C. Fefferman, A. Majda, Geometric 
constraints on potentially singular solutions for the 3-D Euler equations, {\em Commun. in PDE} {\bf 21} (1996), 559-571.
\bibitem{cfbook} P. Constantin, C. Foias, {\em Navier-Stokes equations},
University of Chicago Press, Chicago, 1988.
\bibitem{clm} P. Constantin, P. Lax, A. Majda,  A simple one-dimensional model for the three-dimensional
vorticity equation, {\em Comm. Pure Appl. Math.} {\bf 38} (1985), 715 - 724.
\bibitem{CMT} P. Constantin, A. Majda, and E. Tabak, 
Formation of strong fronts in the 2-D quasi-geostrophic thermal
active scalar, {\it Nonlinearity} {\bf 7} (1994), 1495-1533.
\bibitem{CW1} P. Constantin and J. Wu, Behavior of solutions of 2D
quasi-geostrophic equations, {\it SIAM J. Math. Anal.} {\bf 30}
(1999), 937-948.







\bibitem{Co} D. Cordoba, Nonexistence of simple hyperbolic blow-up
for the quasi-geostrophic equation, {\it Ann. of Math.} {\bf 148} (1998),
1135-1152.

\bibitem{ema} D. Ebin, J. Marsden, Groups of diffeomorphisms and the
motion of an incompressible fluid, {\em Ann. of Math.} (2) {\bf 92} (1970), 102-163.

\bibitem{ey} G. Eyink, Energy dissipation without viscosity in ideal 
hydrodynamics I, Fourier analysis and local transfer, {\em Physica} 
{\bf D 18} (1994), 222-240.

\bibitem{FT} C. Foias and R. Temam, Gevrey class regularity for the solutions
of the Navier-Stokes equations, {\em J. Funct. Anal.}, {\bf 87} (1989),
359-369.

\bibitem{go} J. Gibbon, K. Ohkitani, preprint, to appear in Phys. Fluids.

\bibitem{gft} C. Guillop\'{e}, C. Foias, R. Temam, New a priori estimates for the Navier-Stokes equations in dimension 3, 
{\em Commun. PDE} {\bf 6} (1981),329-359.

\bibitem{mar} D.D. Holm, J. E. Marsden, T. Ratiu, Euler-Poincar\'{e} 
models of ideal fluids with nonlinear dispersion, {\em Phys. Rev. Lett.}
{\bf 349} (1998), 4173-4177.
\bibitem{ka} T. Kato, Nonstationary flows of viscous and ideal fluids in 
R3, {\em J. Funct. Analysis} {\bf 9} (1972), 296-305.

\bibitem{kuz} G. A. Kuzmin, Ideal incompressible hydrodynamics
in terms of the vortex momentum density, {\em Phys. Lett} {\bf 96 A
} (1983), 88-90.

\bibitem{l} J. Leray, Essai sur le mouvement d'un liquide
visqueux emplissant l'espace,  {\em Acta Mathematica} {\bf 63} (1934), 193-248.

\bibitem{pll} P.-L. Lions, {\em Mathematical topics in fluid mechanics}, 
Oxford Lecture Series in Mathematics and its Applications, 3, 
Oxford University Press, New York (1996).

\bibitem{Ma} A. Majda, Vorticity and mathematical theory of 
incompressible fluid flow, {\it Comm. Pure Appl. Math.} {\bf 39}(1986), 
187-220. 


\bibitem{ne} J. Necas, M. Ruzicka, V. Sverak, On Leray's self-similar
solutions of the Navier-Stokes equations, {\em Acta Mathematica} {\bf 176} (1966),
283-294.

\bibitem{o} L. Onsager, Statistical Hydrodynamics, {\em Nuovo Cimento} {\bf 6}(2) (194
9), 279-287.

\bibitem{ose} V. I. Oseledets, On a new way of writing the
  Navier-Stokes equations. The Hamiltonian formalism, {\em Commun. 
Moscow Math. Soc.} (1988), Russ. Math. Surveys {\bf 44} (1989), 210
-211.

\bibitem{R} S. Resnick, {\it Dynamical Problem in Nonlinear Advective Partial Differential Equations}, Ph.D. thesis
University of Chicago , Chicago, 1995.

\bibitem{russo} G. Russo,  P. Smereka, Impulse formulation of the Euler equations:
general properties and numerical methods, {\em J. Fluid Mech.} {\bf  391} (1999), 189-209

\bibitem{serre} D. Serre, La croissance de la vorticite dans les 
ecoulements parfaits incompressibles. {\em C. R. A. S.}, {\bf 328} (1999), 
549-552.

\bibitem{serr} J. Serrin, Mathematical principles of classical fluid 
mechanics, (S. Flugge, C. Truesdell Edtrs.) {\em Handbuch der Physik}, {\bf 8} 
(1959), 125-263.
\bibitem{serrin} J. Serrin, The initial value problem for the Navier-Stokes 
equations, {\em Non-linear Problems}, R. E. Langer edtr., 
Univ. Wisconsin Press, Madison, (1963), 69-98.

\bibitem{jt} J. T. Stuart,  Nonlinear Euler partial differential
equations: singularities in their solution. In {\it Proc. Symp.
 in Honor of C.C. Lin} (D. J. Benney, C. Yuan and F. H. Shu
Edtrs), World Scientific Press, Singapore (1987), 81-95.

\bibitem{t} R. Temam, Some developments on Navier-Stokes Equations in 
the Second Half of the 20th Century, in {\em Developpement des Mathematiques au cours de la seconde moitie du XXeme siecle}, J.P. Pier ed., Birkhauser 
Verlag, Basel 

\bibitem{tbook} R. Temam, {\em Navier-Stokes Equations}, North-Holland, 
Amsterdam-New York, third edition, 1984.

\bibitem{ts} T-P. Tsai, On Leray's self-similar solutions of the Navier-Stokes
equations satisfying local energy estimates, {\em Arch. Rational Mech. Anal.} 
{\bf 143} (1998), 29-51.

\end{thebibliography}
\end{document}